\documentclass[reqno]{amsart}

\usepackage{amsmath,amsthm,amssymb,amsfonts,enumerate,graphicx,color,pstricks}

\numberwithin{equation}{section}

\newcommand{\ti}{\mathrm i}

\begin{document}

\baselineskip 18pt
\larger[2]
\title
[Determinantal elliptic Selberg integrals] 
{Determinantal elliptic Selberg integrals}
\author{Hjalmar Rosengren}
\address
{Department of Mathematical Sciences
\\ Chalmers University of Technology and University of Gothenburg\\SE-412~96 G\"oteborg, Sweden}
\email{hjalmar@chalmers.se}
\urladdr{http://www.math.chalmers.se/{\textasciitilde}hjalmar}

\thanks{Supported by the Swedish Science Research
Council. }

\dedicatory{\Large Dedicated to Christian Krattenthaler on his 60th birthday}

\begin{abstract}
The classical Selberg integral contains a power of the Vandermonde determinant. When that power is a square, it is easy to prove Selberg's identity by interpreting it as a determinant of one-variable integrals. 
We give similar
proofs of summation and transformation formulas for continuous and discrete elliptic Selberg integrals.
In the continuous case, the same proof was previously given by Noumi.
Special cases of the resulting identities have found applications in combinatorics. 
\end{abstract}

\maketitle

\section{Introduction}  
In 1944, Selberg proved the integral evaluation \cite{s2}
\begin{multline}\label{si}\int_{x_1,\dots,x_n=0}^1\, \prod_{1\leq j<k\leq n}|x_j-x_k|^{2c} 
\prod_{j=1}^n x_j^{a-1}(1-x_j)^{b-1}\,dx_j\\ =\prod_{j=1}^n\frac{\Gamma(a+(j-1)c)\Gamma(b+(j-1)c)\Gamma(1+jc)}{\Gamma(a+b+(n+j-2)c)\Gamma(1+c)},
\end{multline}
which had appeared in slightly different form in his earlier paper \cite{s1}.
Here, $\Gamma$ is the classical gamma function, not the elliptic gamma function that will appear below. 
This holds under the natural convergence conditions
$$\operatorname{Re}(a)>0,\quad \operatorname{Re}(b)>0,\quad \operatorname{Re}(c)>-\min\left(\frac 1n,\frac{\operatorname{Re}(a)}{n-1},\frac{ \operatorname{Re}(b)}{n-1}\right). $$
The Selberg integral plays a fundamental role in random matrix theory and analysis on classical groups, and has been generalized in many
directions \cite{fw}.

The general case of \eqref{si} is quite deep. It is instructive to note that in 1963
 Mehta and Dyson \cite{dm} conjectured that
\begin{equation}\label{dmi}\int_{x_1,\dots,x_n=-\infty}^\infty\, \prod_{1\leq j<k\leq n}|x_j-x_k|^{2c} 
\prod_{j=1}^n e^{-x_j^2/2}\,dx_j =(2\pi)^{n/2}\prod_{j=1}^n\frac{\Gamma(1+jc)}{\Gamma(1+c)}.
\end{equation}
Although this was republished as a conjecture several times, no proof was found until it was recognized as a degenerate case of  \eqref{si}  in the late 1970s, see \cite{fw}.

Mehta and Dyson could prove \eqref{dmi} for $c=1/2$, $1$ and $2$, which
are  the most important cases in random matrix theory. Their  discussion of the case $c=1$ 
is just one sentence: ``The case $\beta=2$ is the easiest; one needs only to introduce Hermite polynomials
and exploit their orthogonality properties". For our purpose, it is useful to explain this starting from the algebraic identity \cite{a}
(see \cite{f} for the history of this result)
\begin{multline}\label{ai}\det_{1\leq j,k\leq n}\left(\int f_j(x)g_k(x)\,d\mu(x)\right)\\
=\frac 1{n!}\int \det_{1\leq j,k\leq n}\left(f_j(x_k)\right)
\det_{1\leq j,k\leq n}\left(g_j(x_k)\right)\,d\mu(x_1)\dotsm d\mu(x_n), \end{multline}
which holds for any linear functional
$f\mapsto\int f(x)\,d\mu(x)$. If $f_j$ and $g_j$ are monic polynomials of degree $j-1$, then the determinants on the right are column-equivalent to Vandermonde determinants and we obtain 
\begin{equation}\label{sl}\det_{1\leq j,k\leq n}\left(\int f_j(x)g_k(x)\,d\mu(x)\right)
= \frac 1{n!}\int \prod_{1\leq j<k\leq n}(x_j-x_k)^{2} \,d\mu(x_1)\dotsm d\mu(x_n).\end{equation}
If we now choose $f_j=g_j$ as orthogonal with respect to 
$d\mu$, then the left-hand side of \eqref{sl}  reduces to the product of the squared norms of the first $n$ 
 monic orthogonal polynomials.
The identity \eqref{sl} is then a classical result known to Heine \cite[Cor.~2.1.3]{i}.
The case of Jacobi and Hermite polynomials give the  case $c=1$
of \eqref{si} and \eqref{dmi}, respectively.

There is a less well-known but even more elementary proof of the case $c=1$ of 
\eqref{si}, based on varying the parameters $a$ and $b$. This proof is more relevant to the present work, so we will explain it in detail.
Let $I_{jk}$ denote the one-variable case of \eqref{si}, after
replacing  $(a,b)$ with $(a+j-1,b+n-k)$. By Euler's beta integral evaluation,
$$I_{jk}=\int_0^1 x^{a+j-2}(1-x)^{b+n-k-1}\,dx=
\frac{\Gamma(a+j-1)\Gamma(b+n-k)}{\Gamma(a+b+n+j-k-1)}. $$
Consider the determinant
$D=\det_{1\leq j,k\leq n}(I_{jk})$. It can be identified
with the left-hand side of \eqref{sl}, where   $f_j(x)=x^{j-1}$, $g_j(x)=(1-x)^{n-j}$
and 
$$\int f(x)\,d\mu(x)=\int_0^1 f(x) \,x^{a-1}(1-x)^{b-1}\,dx.$$
Although $g_j$ is not monic of degree $j-1$, the sign changes  coming from
replacing $(x-1)^{j-1}$ by $(1-x)^{n-j}$ cancel, so \eqref{sl} still holds.
Thus, the case $c=1$ of \eqref{si} can be expressed as $n! D$.
This is another instance of the Vandermonde determinant. The gamma functions in the numerator can be pulled out, and the denominator can be expressed as
$$\frac 1{\Gamma(a+b+n+j-k-1)}=\frac{p_{k-1}(j)}{\Gamma(a+b+n+j-2)}, $$
where $p_{k-1}$ is a monic polynomial of degree $k-1$, so that
$$\det_{1\leq j,k\leq n}\left(p_{k-1}(j)\right)=\prod_{1\leq j<k\leq n}(k-j)=\prod_{j=1}^n(j-1)!\,. $$
We conclude  that
$$n!D=\prod_{j=1}^n\frac{\Gamma(a+j-1)\Gamma(b+j-1)j!}{\Gamma(a+b+n+j-2)},$$
which agrees with the right-hand side of \eqref{si} for $c=1$.

The same method works for many variations of the $c=1$ Selberg integral; the  integration may be continuous, infinite discrete or finite, and the integrands may live at the rational, trigonometric or elliptic level.  
(In the most common notation, $c=1$ corresponds to $t=q$ at the trigonometric and elliptic level.)
One can also prove transformation formulas, stating that two Selberg-type integrals are equal. The purpose of the present note is to illustrate this method with two examples:  an elliptic Selberg integral  conjectured by Van Diejen and Spiridonov \cite{ds} and a transformation formula for discrete Selberg integrals
conjectured by Warnaar \cite{w}. Both these identities were first proved by Rains \cite{ra,ri}; the second one also independently by Coskun and Gustafson \cite{cg}.

We are not claiming that our proofs are new, and the present paper should be viewed as expository. 
It is clear from our correspondence with Rains that he is aware of similar proofs. Moreover,
 Noumi \cite{n} gave a determinantal
proof of the  transformation formula stated as \eqref{eit} below. This generalizes the proof of the continuous integration formula given below and is completely parallel to our proof of the discrete transformation formula. 
The main motivation for writing the present note is that we have seen several recent papers where
the case $t=q$ of Warnaar's identities for discrete Selberg integrals are applied \cite{bk,bkw,fkx,ks}
 but the reader is referred to 
work on the general case \cite{cg,ra,rw} for the proof. Even though it 
is known to some experts in the field,
it seems useful to point out to a wider community that much easier proofs exist.   
We also hope that the same method can be used to
find new results. In particular, we think of
quadratic and cubic transformation formulas for $c=1$ Selberg-type integrals, which may
perhaps even admit extensions to general $c$. In this direction, we mention that several quadratic transformations of
  elliptic Selberg integrals are given in \cite{r4,r5}. 
Quadratic summations  for $c=1$ discrete Selberg integrals appear in connection with  tiling problems \cite{c,rsi}.

{\bf Acknowledgement:}  It is a pleasure to dedicate this piece to Christian Krattenthaler, a \emph{virtuoso} of determinants, hypergeometry and much more.
I am very grateful for his patience and support over the years. 
I also thank Masatoshi Noumi, Eric Rains and Ole Warnaar for useful correspondence.

\section{Continuous Selberg integrals}\label{cs}

We recall the standard notation of elliptic hypergeometric functions. We fix two parameters $p$ and $q$ with $|p|,\,|q|<1$, which we suppress from the notation.
 Ruijsenaars' elliptic gamma function \cite{ru}
is given by
$$ \Gamma(z)=\prod_{j,k=1}^n\frac{1-p^{j+1}q^{k+1}/z}{1-p^jq^kz}. $$
It satisfies the functional equation
\begin{equation}\label{gt}\Gamma(qz)=\theta(z)\Gamma(z) \end{equation}
and, more generally,
\begin{equation}\label{gp}\Gamma(q^kz)=(z)_k\,\Gamma(z), \end{equation}
where the theta function and elliptic shifted factorials are given by
$$\theta(z)=\prod_{j=0}^\infty(1-p^j z)\left(1-\frac{p^{j+1}}z\right), \qquad 
(z)_k=\prod_{j=0}^{k-1}\theta(z q^j). $$
Repeated variables in each of these functions is a short-hand for products. For instance,
$$\Gamma(z_1,\dots,z_m)=\Gamma(z_1) \dotsm\Gamma(z_m),$$
$$ \Gamma(z^\pm w^\pm)=\Gamma(zw)\Gamma(z/w)\Gamma(w/z)\Gamma(1/wz).$$
For introductions to elliptic hypergeometric series, we refer to \cite{gr,rs}. We will make
heavy use of  elementary identities that can be found in these sources.

The elliptic Selberg integral is the evaluation
\begin{multline}\label{esi}\frac{C^n}{2^nn!}\int \prod_{1\leq j<k\leq n}\frac{\Gamma(tz_j^{\pm}z_k^{\pm})}{\Gamma(z_j^\pm z_k^\pm)} 
\prod_{j=1}^n\frac{\prod_{k=1}^6\Gamma(t_kz_j^\pm)}{\Gamma(z_j^{\pm 2})}\frac{dz_j}{2\pi\ti z_j}\\
=\prod_{m=1}^n\left(\frac{\Gamma(t^m)}{\Gamma(t)}\prod_{1\leq j<k\leq 6}\Gamma(t^{m-1}t_jt_k)\right),
\end{multline}
where the parameters satisfy the balancing condition $t^{2n-2}t_1t_2t_3t_4t_5t_6=pq$
and
$$C=\prod_{j=1}^\infty(1-p^j)(1-q^j). $$
 If $|t|<1$ and $|t_j|<1$
for all $j$, the integration is  over  $|z_1|=\dots=|z_n|=1$; this condition may be relaxed if
the contour is deformed appropriately. The evaluation
 \eqref{esi} contains the classical Selberg integral \eqref{si} as a limit, see \cite{rl}. 

The case $p=0$ of \eqref{esi} is due to Gustafson \cite{g} and the case $n=1$ to Spiridonov \cite{sp}.
The general case was conjectured by Van Diejen and Spiridonov \cite{ds} and first proved by Rains \cite{ri}.
Another proof follows by combining the results of \cite{ds2,sp2} and a third
 proof is given in   \cite{in}.  For a quantum field theory interpretation of \eqref{esi}, see \cite[\S 12.3.2]{sv}.

The parameter $c$ in \eqref{si} corresponds to $\log_q t$ in \eqref{esi}. In particular, $c=1$ corresponds to $t=q$. 
We proceed to give a simple proof of this special case.
Let $I_{jk}$ denote the case $n=1$ of
\eqref{esi}, after the substitutions 
$(t_1,t_2,t_3,t_4)\mapsto(t_1q^{j-1},t_2q^{n-j},t_3q^{k-1},t_4q^{n-k})$. The balancing condition is then
\begin{equation}\label{bc}q^{2n-2}t_1t_2t_3t_4t_5t_6=pq. \end{equation}
By \eqref{gp}, the resulting integral
can be written as
$$I_{jk}=\frac{C}{2}\int
(t_1z^{\pm})_{j-1}(t_2z^{\pm})_{n-j}(t_3z^{\pm})_{k-1}(t_4z^\pm)_{n-k}
\frac{\prod_{j=1}^6\Gamma(t_jz^{\pm})}{\Gamma(z^{\pm 2})}
\frac{dz}{2\pi\ti z}.$$
Let $D=\det_{1\leq j,k\leq n}(I_{jk})$. Then,    \eqref{ai} gives
$$D=\frac{C^n}{2^nn!}\int d(t_1,t_2)d(t_3,t_4)\prod_{k=1}^n\frac{\prod_{j=1}^6\Gamma(t_jz_k^\pm)}{\Gamma(z_k^{\pm 2})}\frac{dz_k}{2\pi\ti z_k},$$
where
\begin{equation}\label{wdd}d(a,b)=\det_{1\leq j,k\leq n}\big((az_k^\pm)_{j-1}(bz_k^\pm)_{n-j}\big). \end{equation}
By Warnaar's determinant evaluation \cite[Lemma 5.3]{w},
\begin{equation}\label{wd}
d(a,b)=b^{\binom n2}q^{\binom n3}\prod_{j=1}^n(q^{j-n}a/b,q^{n-j}ab)_{j-1}\prod_{1\leq j<k\leq n}z_k^{-1}\theta(z_kz_j^\pm). 
\end{equation}
Note also that
$$  \prod_{1\leq j<k\leq n}
\big(z_k^{-1}\theta(z_kz_j^\pm)\big)^2=\prod_{1\leq j<k\leq n}\frac{\Gamma(qz_j^{\pm}z_k^{\pm})}{\Gamma(z_j^\pm z_k^\pm)}. $$
This gives
\begin{align}\nonumber D&=(t_2t_4)^{\binom n2}q^{2\binom n3}\prod_{j=1}^n
(q^{j-n}t_1/t_2,q^{n-j}t_1t_2,q^{j-n}t_3/t_4,q^{n-j}t_3t_4)\\
 \label{di}&\quad\times \frac{C^n}{2^n n!}\int \prod_{1\leq j<k\leq n}\frac{\Gamma(qz_j^{\pm}z_k^{\pm})}{\Gamma(z_j^\pm z_k^\pm)} 
\prod_{j=1}^n\frac{\prod_{k=1}^6\Gamma(t_kz_j^\pm)}{\Gamma(z_j^{\pm 2})}\frac{dz_j}{2\pi\ti z_j},\end{align}
where we recognize the integral as the case $t=q$ of \eqref{esi}.

On the other hand, the case $n=1$ of \eqref{esi} (that is, Spiridonov's elliptic beta integral) gives
\begin{align*}I_{jk}&=\Gamma(t_1t_2q^{n-1},t_1t_3q^{j+k-2},t_1t_4q^{n+j-k-1},t_1t_5q^{j-1},t_1t_6q^{j-1})\\
&\quad\times\Gamma(t_2t_3q^{n-j+k-1},t_2t_4q^{2n-j-k},t_2t_5q^{n-j},t_2t_6q^{n-j},t_3t_4q^{n-1})\\
&\quad\times\Gamma(t_3t_5q^{k-1},t_3t_6q^{k-1},t_4t_5q^{n-k},t_4t_6q^{n-k},t_5t_6).
\end{align*}
Most of the factors are independent of either $j$ or $k$ and can thus be pulled out of the determinant.
Using again \eqref{gp}, we are left with
\begin{align*}
D&=\Gamma(t_1t_2q^{n-1},t_3t_4q^{n-1},t_5t_6)^n\prod_{m=1}^n\prod_{\substack{1\leq j<k\leq 6\\
(j,k)\neq (1,2),\,(3,4),\,(5,6)}}\Gamma(t_jt_kq^{m-1})\\
&\quad\times\det_{1\leq j,k\leq n}\left((t_1t_3q^{k-1},t_1t_4q^{n-k})_{j-1}(t_2t_3q^{k-1},t_2t_4q^{n-k})_{n-j}\right).
\end{align*}
The final determinant is of the form \eqref{wdd}, with $a=t_1\sqrt{t_3t_4q^{n-1}}$,
$b=t_2\sqrt{t_3t_4q^{n-1}}$, $z_k=q^{k-1}\sqrt{q^{1-n}t_3/t_4}$. Using \eqref{wd} and also
 \eqref{bc} to write $(q^{n-j}ab)_{j-1}=(t_5t_6)_{j-1}$,
we obtain after simplification
\begin{align*}
D&=(t_2t_4)^{\binom n2}q^{2\binom n3}\Gamma(t_1t_2q^{n-1},t_3t_4q^{n-1})^n\prod_{m=1}^n\prod_{\substack{1\leq j<k\leq 6\\
(j,k)\neq (1,2),\,(3,4)}}\Gamma(t_jt_kq^{m-1})\\
&\quad\times\prod_{j=1}^n(q,q^{j-n}t_1/t_2,q^{j-n}t_3/t_4)_{j-1}.
\end{align*}
Comparing this with \eqref{di} yields the case $t=q$ of \eqref{esi}.

Essentially the same proof works for  the case $t=q$ of
Rains' integral transformation \cite{ri}
\begin{multline}\label{eit}\int \prod_{1\leq j<k\leq n}\frac{\Gamma(tz_j^{\pm}z_k^{\pm})}{\Gamma(z_j^\pm z_k^\pm)} 
\prod_{j=1}^n\frac{\prod_{k=1}^4\Gamma(t_kz_j^\pm,u_kz_j^{\pm})}{\Gamma(z_j^{\pm 2})}\frac{dz_j}{2\pi\ti z_j}\\
\begin{split}&=\prod_{m=1}^n\prod_{1\leq j<k\leq 4}\Gamma(t^{m-1}t_jt_k,t^{m-1}u_ju_k)\\
&\quad\times\int \prod_{1\leq j<k\leq n}\frac{\Gamma(tz_j^{\pm}z_k^{\pm})}{\Gamma(z_j^\pm z_k^\pm)} 
\prod_{j=1}^n\frac{\prod_{k=1}^4\Gamma(t_kvz_j^\pm,u_kv^{-1}z_j^{\pm})}{\Gamma(z_j^{\pm 2})}\frac{dz_j}{2\pi\ti z_j},
\end{split}\end{multline}
where $v^2=pq/t^{n-1}t_1t_2t_3t_4=t^{n-1}u_1u_2u_3u_4/pq$. 
The details can be found in \cite{n}, where the integral \eqref{eit} is interpreted as a tau function for the elliptic Painlev\'e equation.

\section{Discrete Selberg integrals}

The integral evaluation \eqref{esi} and transformation \eqref{eit} have 
analogues for finite sums, which were conjectured by Warnaar \cite{w} prior to the discovery of the continuous
versions. Warnaar's summation can be obtained from the integral evaluation \eqref{esi}
through residue calculus \cite{ds} (presumably, a similar argument applies to the transformations).
The conjectured summation was proved in \cite{rw}, see also \cite{ins}, and the more general transformation 
in \cite{cg,ra}. 

As one would expect, the case $t=q$ of Warnaar's identities admit simple determinantal proofs. We will focus on the transformation, which can be written as~\cite{sc}
\begin{multline}\label{mbt}
\sum_{0\leq x_1<x_2<\dots<x_n\leq N}\,\prod_{1\leq j<k\leq n}\big(q^{x_j}\theta(q^{x_k-x_j})\theta(aq^{x_j+x_k})\big)^2\\
\begin{split}&\quad\times\prod_{j=1}^n\left(\frac{\theta(aq^{2x_j})}{\theta(a)}
\frac{(a,b,c,d,e,f,g,q^{-N})_{x_j}}{(q,aq/b,aq/c,aq/d,aq/e,aq/f,aq/g,aq^{N+1})_{x_j}}
\,q^{x_j}\right)\\
&=\left(\frac a\lambda\right)^{(N+1-n)n}\frac{(aq)_N^n}{(\lambda q)_N^n}
\prod_{j=1}^n\frac{(b,c,d)_{j-1}(\lambda q/e,\lambda q/f,\lambda q/g)_{N+1-j}}{(\lambda b/a,\lambda c/a,\lambda d/a)_{j-1}(aq/e,aq/f,aq/g)_{N+1-j}}\\
&\quad\times\sum_{0\leq x_1<x_2<\dots<x_n\leq N}\,\prod_{1\leq j<k\leq n}\big(q^{x_j}\theta(q^{x_k-x_j})\theta(\lambda q^{x_j+x_k})\big)^2\\
&\quad\times\prod_{j=1}^n\left(\frac{\theta(\lambda q^
{2x_j})}{\theta(\lambda)}
\frac{(\lambda,\lambda b/a,\lambda c/a,\lambda d/a,e,f,g,q^{-N})_{x_j}}{(q,aq/b,aq/c,aq/d,\lambda q/e,\lambda q/f,\lambda q/g,\lambda q^{N+1})_{x_j}}
\,q^{x_j}\right),
\end{split}\end{multline}
where $bcdefg=q^{4+N-2n}a^3$ and $\lambda=a^2q^{2-n}/bcd$. 
When 
$aq=cd$,  the factor $(\lambda b/a)_{x_n}=(q^{1-n})_{x_n}$ on the right-hand side vanishes unless $x_n\leq n-1$,
so the sum reduces to the term with $(x_1,\dots,x_n)=(0,1,\dots,n-1)$. 
After a change of variables,
this gives the case $t=q$ of Warnaar's discrete Selberg integral, namely,
\begin{multline}\label{mbs}
\sum_{0\leq x_1<x_2<\dots<x_n\leq N}\,\prod_{1\leq j<k\leq n}\big(q^{x_j}\theta(q^{x_k-x_j})\theta(aq^{x_j+x_k})\big)^2\\
\begin{split}&\quad\times\prod_{j=1}^n\left(\frac{\theta(aq^{2x_j})}{\theta(a)}
\frac{(a,b,c,d,e,q^{-N})_{x_j}}{(q,aq/b,aq/c,aq/d,aq/e,aq^{N+1})_{x_j}}
\,q^{x_j}\right)\\
&=b^{n(N+1-n)}q^{\frac 13\,n(n-1)(3N+1-2n)}(aq)_N^n\\
&\quad\times\prod_{j=1}^n\frac{(q,b,c,d,e,q^{-N})_{j-1}(aq^{2-j}/bc,aq^{2-j}/bd,aq^{2-j}/be)_{N+1-n}}{(aq/b,aq/c,aq/d,aq/e)_{N+1-j}},
\end{split}\end{multline}
which holds for $bcde=q^{N+3-2n}a^2$.

We will give a simple  proof of \eqref{mbt}, which  is completely parallel to the continuous case. 
We need the case $n=1$, which is the one-variable elliptic Bailey transformation. It first appeared (rather implicitly and with some restrictions on the parameters) in the work of Date et al.\ on Baxter's elliptic solid-on-solid model \cite{d}
and was proved in general by Frenkel and Turaev \cite{ft}, see \cite{gr,rs} for more elementary proofs.

Let $S_{jk}$ denote 
the case $n=1$ of \eqref{mbt}, after the substitutions
$(b,c,e,f)\mapsto(bq^{j-1},cq^{n-j},eq^{k-1},fq^{n-k})$. After some elementary manipulation, we may write 
\begin{align*}
S_{jk}&=\sum_{x=0}^N\frac{\theta(aq^{2x})}{\theta(a)}
\frac{(a,b,c,d,e,f,g,q^{-N})_{x}}{(q,aq/b,aq/c,aq/d,aq/e,aq/f,aq/g,aq^{N+1})_{x}}\,q^{x(2n-1)}\\
&\quad\times\frac{(bq^x,bq^{-x}/a)_{j-1}(cq^x,cq^{-x}/a)_{n-j}(eq^x,eq^{-x}/a)_{k-1}(fq^x,fq^{-x}/a)_{n-k}}
{(b,b/a)_{j-1}(c,c/a)_{n-j}(e,e/a)_{k-1}(f,f/a)_{n-k}}.
\end{align*}
Let $D=\det_{1\leq j,k\leq n}(S_{jk})$. 
We expand $D$ using the Cauchy--Binet identity
$$\det_{1\leq j,k\leq n}\left(\sum_{x=0}^N a_{jx}b_{kx}\right)
=\sum_{0\leq x_1<x_2<\dots< x_n\leq N}\det_{1\leq j,k\leq n}(a_{j,x_k})\det_{1\leq j,k\leq n}(b_{j,x_k}).
 $$
This is a special case of \eqref{ai}, where symmetry is used to restrict the range of summation.
It follows that
\begin{multline*}
D=\frac 1{\prod_{j=1}^n(b,b/a,c,c/a,e,e/a,f,f/a)_{j-1}}
\\
\times\sum_{0\leq x_1<x_2<\dots< x_n\leq N}\,\prod_{k=1}^n\left(\frac{\theta(aq^{2x_k})}{\theta(a)}\frac{(a,b,c,d,e,f,g,q^{-N})_{x_k}\,q^{x_k(2n-1)}}{(q,aq/b,aq/c,aq/d,aq/e,aq/f,aq/g,aq^{N+1})_{x_k}}\right)\\
\times \delta(b,c)\delta(e,f),\end{multline*}
where
$$\delta(b,c)= \det_{1\leq j,k\leq n}\left((bq^{x_k},bq^{-x_k}/a)_{j-1}(cq^{x_k},cq^{-x_k}/a)_{n-j}\right).$$
This determinant  is of the form \eqref{wdd}, with $(a,b,z_k)$ replaced by $(b/\sqrt a,c/\sqrt{a},\sqrt a q^{x_k})$. 
Using \eqref{wd} and simplifying, we find that $D$ equals the left-hand side of \eqref{mbt} times
\begin{equation}\label{pf}\left(\frac{cf}{a^2}\right)^{\binom n2}q^{2\binom n3}\prod_{j=1}^n
\frac{(q^{j-n}b/c,q^{n-j}bc/a,q^{j-n}e/f,q^{n-j}ef/a)_{j-1}}{(b,b/a,c,c/a,e,e/a,f,f/a)_{j-1}}.\end{equation}

Repeating the same computation but starting from the alternative expression
\begin{align*}
S_{jk}&=\left(\frac a\lambda\right)^{N}
\frac{(aq,\lambda q^{2-k}/e,\lambda q^{1-n+k}/f,\lambda q/g)_{N}}{(\lambda q,aq^{2-k}/e,aq^{1-n+k}/f,aq/g)_{N}}\\
&\quad\times\sum_{x=0}^N\frac{\theta(\lambda q^{2x})}{\theta(\lambda)}
\frac{(\lambda,\lambda b/a,\lambda c/a,\lambda d/a,e,f,g,q^{-N})_{x}}{(q,aq/b,aq/c,aq/d,\lambda q/e,\lambda q/f,\lambda q/g,\lambda q^{N+1})_{x}}\,q^{x(2n-1)}\\
&\quad\times\frac{(\lambda bq^x/a,bq^{-x}/a)_{j-1}(\lambda cq^x/a,cq^{-x}/a)_{n-j}(eq^x,eq^{-x}/\lambda)_{k-1}(fq^x,fq^{-x}/\lambda)_{n-k}}
{(\lambda b/a,b/a)_{j-1}(\lambda c/a,c/a)_{n-j}(e,e/\lambda)_{k-1}(f,f/\lambda)_{n-k}},
\end{align*}
we obtain after simplification the same prefactor \eqref{pf} times the right-hand side of \eqref{mbt}.
This completes the proof of \eqref{mbt}.


\begin{thebibliography}{99}
\bibitem[A]{a} C.\ Andr\'eief, \emph{Note sur une relation entre les int\'egrales d\'efinies des produits des fonctions},
M\'em.\ Soc.\ Sci.\ Phys.\ Nat.\ Bordeaux 2 (1886), 1--14.
\bibitem[BK]{bk} G.\ Bhatnagar and C.\ Krattenthaler,
\emph{The determinant of an elliptic Sylvesteresque matrix}, SIGMA 14 (2018), 052.
\bibitem[BKW]{bkw} R.\ P.\ Brent, C.\ Krattenthaler and S.\ O.\ Warnaar,
\emph{Discrete analogues of Macdonald--Mehta integrals}, J.~Combin.\ Theory Ser.\ A  144  (2016), 80--138. 
\bibitem[C]{c} M.\ Ciucu, T.\ Eisenk\"olbl, C.\ Krattenthaler and D.\ Zare, \emph{Enumeration of lozenge tilings of hexagons with a central triangular hole}, J.\ Combin.\ Theory Ser. A  95  (2001),  251--334. 
\bibitem[CG]{cg}
H.\ Coskun and R.\ A.\ Gustafson, \emph{Well-poised Macdonald functions $W_\lambda$ and Jackson coefficients $\omega_\lambda$ on $BC_n$}, in  Jack, Hall--Littlewood and Macdonald polynomials, Amer.\ Math.\ Soc., 2006,  
127--155.
\bibitem[D]{d} E.\ Date, M.\ Jimbo, A.\ Kuniba, T.\ Miwa, and M.\ Okado, \emph{Exactly solvable SOS models.\ II.\ Proof of the star-triangle relation and combinatorial identities}, in Conformal Field Theory and Solvable Lattice Models, Academic Press, 1988, 17--122.
\bibitem[DS1]{ds} J.\ F.\ van Diejen and V.\ P.\ Spiridonov, 
\emph{An elliptic Macdonald--Morris conjecture and multiple modular hypergeometric sums},
Math.\ Res.\ Lett.\  7  (2000),  729--746. 
\bibitem[DS2]{ds2} J.\ F.\ van Diejen and V.\ P.\ Spiridonov, 
\emph{Elliptic Selberg integrals}, Internat.\ Math.\ Res.\ Notices 2001, 1083--1110. 
\bibitem[FKX]{fkx} H.\ Feng, C.\ Krattenthaler and Y.\ Xu,
\emph{Best polynomial approximation on the triangle}, arXiv:1711.04756.
\bibitem[F]{f} P.\ J.\ Forrester , \emph{Meet Andr\'eief, Bordeaux 1886, and Andreev, Kharkov 1882--83},
arXiv:1806.10411.
\bibitem[FW]{fw} 
P.\ J.\ Forrester and S.\ O.\  Warnaar, \emph{The importance of the Selberg integral}, Bull.\ Amer.\ Math.\ Soc.\  45  (2008),  489--534.
\bibitem[FT]{ft} I.\ B.\ Frenkel and V.\ G.\ Turaev, \emph{Elliptic solutions of the Yang--Baxter equation and modular hypergeometric functions}, in  The Arnold--Gelfand Mathematical Seminars, Birkh\"auser, 1997, 171--204. 
\bibitem[GR]{gr} G.\ Gasper and M.\ Rahman, Basic Hypergeometric Series, Second edition, Cambridge University Press, 2004.
\bibitem[G]{g} R.\ A.\ Gustafson,
\emph{Some $q$-beta integrals on $\mathrm{SU}(n)$  and $\mathrm{Sp}(n)$  that generalize the Askey--Wilson and Nasrallah--Rahman integrals}, 
SIAM J.\ Math.\ Anal.\  25  (1994),  441--449.  
\bibitem[I]{i} M.\ E.\ H.\ Ismail, Classical and Quantum Orthogonal Polynomials in One Variable, 
 Cambridge University Press,  2005.
\bibitem[IN1]{ins} M.\ Ito and M.\ Noumi, \emph{Derivation of a $BC_n$ elliptic summation formula via the fundamental invariants}, Constr.\ Approx.\ 45 (2017), 33--46.
\bibitem[IN2]{in} M.\ Ito and M.\ Noumi, \emph{Evaluation of the $BC_n$ elliptic Selberg integrals via the fundamental invariants}, Proc.\ Amer.\ Math.\ Soc.\ 145 (2017), 689--703.
\bibitem[KS]{ks} C.\ Krattenthaler and M.\ Schlosser,
\emph{The major index generating function of standard Young tableaux of shapes of the form ``staircase minus triangle"}, in Ramanujan 125,  Amer.\ Math.\ Soc., 2014, 111--122.
\bibitem[MD]{dm} M.\ L.\ Mehta and F.\ J.\ Dyson, 
\emph{Statistical theory of the energy levels of complex systems.\ V.} 
J.\ Math.\ Phys.~ 4  (1963), 713--719. 
\bibitem[N]{n} M.\ Noumi, \emph{Remarks on $\tau$-functions for the difference Painlev\'e equations of type $E_8$},
arXiv:1604.06869.
\bibitem[R1]{ra} E.\ M.\ Rains, \emph{$BC_n$-symmetric abelian functions}, Duke Math.\ J.\ 135 (2006), 99--180.
\bibitem[R2]{rl} E.\ M.\ Rains, \emph{Limits of elliptic hypergeometric integrals}, Ramanujan J.\ 18 (2009), 257--306.
\bibitem[R3]{ri} E.\ M.\ Rains, \emph{Transformations of elliptic hypergeometric integrals}, Ann.\ Math.\ 171 (2010), 169--243.
\bibitem[R4]{r4} E.\ M.\ Rains, \emph{Elliptic Littlewood identities}, J.\ Combin.\ Theory A 119 (2012), 1558--1609.
\bibitem[R5]{r5} E.\ M.\ Rains, \emph{Multivariate quadratic transformations and the interpolation kernel},
SIGMA  14  (2018), 019. 
\bibitem[Ro1]{rw} H.\ Rosengren, \emph{A proof of a multivariable elliptic summation formula conjectured by Warnaar}, in q-Series with Applications to Combinatorics, Number Theory, and Physics, Amer.\ Math.\ Soc., 2001, 193--202. 
\bibitem[Ro2]{rsi} H.\ Rosengren, \emph{Selberg integrals, Askey--Wilson polynomials and lozenge tilings of a hexagon with a triangular hole}, J.\ Combin.\ Theory A 138 (2016), 29--59.
\bibitem[Ro3]{rs} H.\ Rosengren, \emph{Elliptic hypergeometric functions}, arXiv:1608.06161.
\bibitem[Ru]{ru} S.\ N.\ M.\ Ruijsenaars,
\emph{First order analytic difference equations and integrable quantum systems},
J.\ Math.\ Phys.\  38  (1997),   1069--1146. 
\bibitem[S]{sc} M.\ Schlosser, \emph{Elliptic enumeration of nonintersecting lattice paths}, J.\ Combin.\ Theory A 114 (2007), 505--521. 
\bibitem[S1]{s1} A.\ Selberg, \emph{\"Uber einen Satz von A.\ Gelfond}, 
Arch.\ Math.\ Naturvid.\  44,  (1941), 159--170. 
\bibitem[S2]{s2} A.\ Selberg, \emph{Bemerkninger om et multipelt integral}, Norsk Mat.\ Tidsskr.\ 26 (1944), 71--78.
\bibitem[Sp1]{sp} V.\ P.\ Spiridonov, \emph{On the elliptic beta function}, Russian Math.\ Surveys 56 (2001), 185--186. 
\bibitem[Sp2]{sp2} V.\ P.\ Spiridonov, \emph{Short proofs of the elliptic beta integrals},
Ramanujan J.\ 13 (2007), 265--283. 
 \bibitem[SV]{sv} V.\ P.\ Spiridonov and G.\ S.\ Vartanov, \emph{Elliptic hypergeometry of supersymmetric dualities}, Comm.\ Math.\ Phys.\ 304 (2011), 797--874. 
\bibitem[W]{w} S.\ O.\ Warnaar, \emph{Summation and transformation formulas for elliptic hypergeometric series}, Constr.\ Approx.\ 18 (2002), 479--502. 
\vskip 3mm
\end{thebibliography}
 \end{document}